\documentclass[12pt]{article}
\usepackage{latexsym,amsmath,amssymb}
\usepackage{amscd,amsmath,amsfonts,graphicx}
\usepackage{amsthm}

\newcommand{\R}{\mathbb{R}}

\newcommand{\E}{{\mathsf E}}
\newcommand{\N}{\mathbb{N}}
\newcommand{\e}{\hbox{e}}

\begin{document}

\begin{center}
{\Large {\bf New Sufficient Conditions for Moment Determinacy\\
via Probability Density Tails}}\footnote{The first submission,
see arXiv 2507.05615v1 [mathPR] 8 Jul 2025, was slightly improved, expanded and reformatted,
and it was published in: {\it Mathematics} (MDPI) {\bf 2025}, 13, 2671. It is available at:  http://doi.org/10.3390/math13162671}
\end{center}

\vspace{0.4cm}
\centerline{{\sc {Gwo Dong Lin}\,$^{a}$ \ and \ Jordan M. Stoyanov\,$^{b}$}}

\vspace{0.3cm}{\small
$^{a}$ Institute of Statistical Science, Academia Sinica, Taipei 11529, Taiwan, ROC }

\hspace{0.20cm} {\small E-mail: \ gdlin@stat.sinica.edu.tw }

$^{b}$  {\small Institute of Mathematics \& Informatics, Bulgarian Academy of Sciences, 1113 Sofia, Bulgaria}
	
\hspace{0.17cm}  {\small E-mail: \ stoyanovj@gmail.com}

\vspace{0.3cm}\noindent
{\small	{\bf Abstract}. One of the ways to characterize a probability distribution is to show that it is moment-determinate, uniquely determined by knowing all its  moments.
The uniqueness, in the absolutely continuous case, depends entirely on the behaviour of the tails of the  density function $f$. We find and exploit  a condition, (D), in terms only of $f$ which is of a `general' form and easy to check. Condition (D), showing the `speed' for $f$ to tend to zero, is
sufficient to conclude the moment determinacy.
We establish a series of theorems and corollaries in both Stieltjes and Hamburger cases and provide an interesting illustrative example. The results in this paper are either new or extend some recently published results.

\vspace{0.1cm}\noindent

\vspace{0.1cm}\noindent
}

\section{Introduction}

In this paper, we address one of the aspects, namely, the uniqueness, in the classical moment problem
 ([1], [2] and [3]). Thus, the interest is in
 conditions under which  a measure, in particular, a probability distribution,
is characterized uniquely by the sequence of all its moments.

We assume that $(\Omega, \mathcal{F}, \mathsf{P})$ is an underlying probability space on which are
defined all random variables considered in this paper. If $X$ is a random variable with distribution $F$,
we write $X \sim F$ and deal with the two possible cases:

\begin{enumerate}
\item[(i)] \ (Hamburger case) $X$ takes values in the real line $\R=(-\infty, \infty).$
\item[(ii)] \ (Stieltjes case) $X$ takes values in the half-real line $\R_+=[0, \infty)$.
\end{enumerate}

For each range of values, $\R$ or $\R_+$, let us assume that $X\sim F$
has finite moments of any positive integer order. This means that any power $|X|^n$ of $X$   is $\mathsf{P}$-integrable for
all $n\in\N :=\{1, 2,\ldots\},$ in which case we have finite moments and the moment sequence, denoted, respectively, as follows:
\[
m_n:={\E}[X^n], \ n\in\N, \quad  \mbox{ and } \quad \{m_n\}_{n=1}^{\infty}.
\]

For any random variable $X \sim F$ with finite moments, there are two possibilities: either $F$ is uniquely determined by its moments,
and we say that $F$ is M-determinate (M-det), or $F$ is non-unique, M-indeterminate (M-indet).
In the latter case, there is at least one distribution $G \ne F$ with the same moments
as $F$.
These notions/properties are equally used in both cases, Hamburger and Stieltjes.

It is clear that if a distribution is M-det on $\R,$ then it is also M-det on $\R_+$. However, the converse is not true in general. It is possible that a distribution $F$ is M-det on $\R_+$ (Stieltjes case); however, it is M-indet on $\R$ (Hamburger case).
 The  meaning of these is that there is no
distribution $G$ on $\R_+$ such that $G \ne F$ and has the same moment sequence, but there does exist another
distribution, say $\tilde G$  on $\R$, \ ${\tilde G} \ne F$ which has the same moment sequence as $F.$
This may happen only for some  discrete distributions; for details, see, e.g., [4].

Most famous and useful conditions which guarantee the M-determinacy were found 100 years ago in [5]. Here are the two statements:

\begin{enumerate}
\item[(i)]	 (Hamburger case) \quad For $X\sim F$ on $\R,$ 
\[
({\rm C_H})~~~~~\sum_{n=1}^{\infty}(\E[X^{2n}])^{-1/(2n)}=\sum_{n=1}^{\infty}m_{2n}^{-1/(2n)}=\infty \quad \Rightarrow \quad X \ \mbox{ is M-det  on } \ \R.
\]
\item[(ii)]	 (Stieltjes case) \quad For  $X\sim F$ on $\R_+,$ 
\[ ({\rm C_S})~~~~\sum_{n=1}^{\infty}(\E[X^{n}])^{-1/(2n)}=\sum_{n=1}^{\infty}m_{n}^{-1/(2n)}=\infty   \quad \Rightarrow \quad X \ \mbox{ is M-det  on } \ \R_+.
\]
\end{enumerate}

Traditionally, $({\rm C_H})$ and $({\rm C_S})$
are each called Carleman's condition and it is well-known that, in a sense, this condition is the `best' sufficient condition for M-det;
see Section 11 in  [6] or [7]. Notice that Carleman's condition is in the group of the
so-called `checkable' sufficient conditions for M-determinacy. It is useful to mention that
several `checkable' sufficient or `checkable'  necessary conditions for either M-det
or M-indet  can be found in the recent works  [8] and [9].

{In this paper, we focus on distributions which are arbitrary in a neighbourhood of  zero, say, for $x \in (-x_0, x_0)$ for some $x_0>0$; however, they are absolutely continuous on a subset of $\R$ outside that interval, so for $X \sim F$, there is a density $f = F'$. For the M-det property, it is important to note the behaviour of  the `probability density tail(s)', $f(x)$ as $|x| \to \infty.$ 
We have two tails on $\R$ and one tail on $\R_+$.  Since $f$ is a density, the rate of its decreasing to zero is related to the rate of growing to infinity of the moments $m_n$ as $n \to \infty,$ which in turn is decisive for the divergence of the Carleman's series; see $\rm (C_H), \ (C_S)$ above.

  Recently, an interesting condition was introduced in [10].
   The idea, expressed in our words,  is to assume that  $f(x) > 0$ on $\R$, use a positive function $\phi$,  and for large $|x|$, compare the `small values'   $f(x \pm \phi(|x|))$ and $f(x)$. This is carried out in terms of
  the density `dropping speed', as we call the behaviour of the
  ratio $f(x \pm \phi(|x|))/f(x)$ as $x \to \pm \infty$.

Specifically, these authors  showed that if  for some constants $a>0, \alpha\in[0,1],$
\ with ${\rm{sign}}(x)=1$ if $x \geq 0$ or $-1$ if $x<0$, the density $f$ satisfies the inequality
\begin{eqnarray*}
\limsup_{|x|\to\infty}\frac{f(x+a\,{\rm{sign}}(x)(\log |x|)^{\alpha})}{f(x)}<1,
\end{eqnarray*}
 then
$X\sim F$ satisfies Carleman's condition $({\rm C_H})$ and is M-det. However, in their proof, the authors assume implicitly that
the underlying distribution $F$ is symmetric about zero and
is absolutely continuous on the whole real line $\R$  (see, e.g., relation (3.8) in the proof of Theorem 2.1 in [10]).

In the present paper, we extend and/or slightly modify the findings of [10].
We suggest considering  `more general' functions $\phi$ and find precise conditions under which one single asymptotic property of the density $f$ implies M-det. Symbolically:

\[
{\bf Condition \ (D)} \hspace{1.0cm} \limsup_{|x|\to\infty}\frac{f(x+{\rm{sign}}(x)\phi(|x|))}{f(x)}<1 \quad \Rightarrow
\quad X \sim F \ \mbox{ is M-det.}
\]

The structure of the paper is as follows. In Section \ref{sec2}, we state the main results, three theorems, and four corollaries. The needed lemmas for the proofs are given in Section \ref{sec3}.
The complete proofs of the main results are provided in Section \ref{sec4}. Comments and
comprehensive illustrative examples are given in Section \ref{sec5}.

\section{Main Theorems and Corollaries}\label{sec2}

{\bf Assumption 1.}
We start with a non-negative and differentiable  function  $\phi$ on $[x_0,\infty)$ for some $x_0> 0.$
Suppose that the function $y = y(x) =x+\phi(x)$ is strictly increasing in $x$ and let $x(y)$ be its inverse function on $[y_0,\infty)$ with $y_0=y(x_0).$
Define the function $\varphi(y)=y-x(y)$ and assume that it satisfies the following conditions, (a) -- (c):

(a) \ $\varphi^{\prime}(y)\in[0,1]$ for $y\ge y_*$;

(b) \ $\varphi(y)\le C_+\log y$ for $y\ge y_*$;

(c) \ $y\varphi^{\prime}(y)\le C_+$ for $y\ge y_*$,
where  $C_+>0$ and $y_*\ge y_0$ are two constants.

\smallskip
We are going to use below the following two notations:
\[
\R_{x_0}:=\R\setminus (-x_0, x_0) \quad \mbox{ and } \quad \R_{+, x_0}:=\R_+\setminus [0,  x_0) \ \mbox{ for some } \ x_0 >0.
\]

\subsection{Statement and a Corollary in the Hamburger Case}

{\bf Theorem 1.}
{\it Let $X\sim F$ on $\R$ have finite moments  and  $F$ be absolutely continuous
on the set $\R_{x_0}$ with $F^{\prime}(x)=f(x)>0$ on $\R_{x_0}.$
Under these assumptions and Assumption 1 for the functions $\phi$ and $\varphi$, assume further that the following relation is satisfied:
\begin{eqnarray}\gamma_1(f,\phi):=\limsup_{|x|\to\infty}\frac{f(x+{\rm{sign}}(x)\phi(|x|))}{f(x)}<1.
\end{eqnarray}
Then the moments of $X$ satisfy Carleman's condition $({\rm C_H})$ and hence $X$ is M-det on $\R.$ Moreover, both $X^2$ and $|X|$ are also M-det on $\R$ and on $\R_+.$}

\vspace{0.1cm}
Notice that in Theorem 1, $X$ is on $\R$, and the first claim is about M-det on $\R$.
The two other random quantities, $X^2$ and $|X|$, even being non-negative, can be considered also on $\R$,
not only on $\R_+.$ See again the comment in the Introduction with  a reference to [4].
Hence there is no reason to worry if there is a claim for a non-negative random variable and a property on the whole real line $\R$.

\smallskip
Consider now specific functions $\phi$ as described in Assumption 1 and derive explicit corollaries. Let us introduce the following three functions, each being a possible choice for $\phi (x)$ with constants $a>0$ and $\alpha \in [0,1]:$
\begin{eqnarray}
\phi_1 (x) &=& a (\log x)^{\alpha}, \quad  x >1; \\
\phi_2 (x) &=& a(\log x)^{\alpha}+\log\log x, \quad x>\e; \\
\phi_3 (x) &=& a(\log x)^{\alpha}/\log\log x, \ x > \e.
\end{eqnarray}

\vspace{0.1cm}\noindent
{\bf Corollary 1.}
{\it Let $X\sim F$ on $\R$ have all finite moments and  $F$ be absolutely continuous
on $\R_{x_0}$ for some $x_0> \e$ with $F^{\prime}(x)=f(x)>0$ on $\R_{x_0}.$ Let $a>0, \ \alpha\in[0,1].$ 
Suppose that for any choice of $\phi (x)$ from $(2)-(4)$, the following relation is valid:
\vspace{-0.2cm}
\[
\limsup_{|x|\to\infty}\frac{f(x+{\rm{sign}}(x)\phi (|x|))}{f(x)} < 1.
\]
Then, the moments of $X$ satisfy Carleman's condition $({\rm C_H})$ and hence $X$ is M-det. Moreover, $X^2$ and $|X|$ are also M-det on $\R$ and on $\R_+$. }

\smallskip
Let us mention that  $\phi (x) = \phi_1 (x)$ is the only function used in
[10].  Here we suggest a slight extension, e.g., adding to $\phi_1 (x)$, or dividing it by, 
another `slowly varying function'.    Notice also that, for large $x$,  the following relation holds:

\smallskip
\centerline{$\phi_3(x) < \phi_1 (x) < \phi_2 (x).$ }

\subsection{Statements and Corollaries in the Stieltjes Case}

{\bf Theorem 2.}
{\it Let $Y\sim G$ on $\R_+$ have all finite moments and  $G$ be absolutely continuous
on $\R_{+,x_0}$ with $G^{\prime}(x)=g(x)>0$ on $\R_{+,x_0}.$ Under these conditions
and Assumption 1 for the functions $\phi$ and $\varphi$, let us further assume that
\begin{eqnarray}\lim_{x\to\infty}\frac{\phi(x)}{x}=0, \quad \gamma_2(g,\phi):=\limsup_{x\to\infty}\frac{g((x+\phi(x))^2)}{g(x^2)}<1.
\end{eqnarray}
Then, the moments of $Y$ satisfy Carleman's condition $({\rm C_S})$ and hence $Y$ is M-det on $\R_+.$}

\vspace{0.2cm}\noindent
{\bf Corollary 2.}
{\it Let $Y\sim G$ on $\R_+$ have all finite moments  and  $G$ be absolutely continuous
on $\R_{+,x_0}$  with $x_0>\e$, $G^{\prime}(x)=g(x)>0$ on $\R_{+,x_0}.$ Suppose that $\phi (x)$ is any one of the functions 
$(2)-(4)$,  with $a>0,\ \alpha\in[0,1]$, and that the following relation is satisfied:
\begin{eqnarray*}
	\limsup_{x\to\infty}\frac{g((x+ \phi (x))^2)}{g(x^2)}<1.
\end{eqnarray*}
Then, the moments of $Y$ satisfy Carleman's condition $({\rm C_S})$ and hence $Y$ is M-det on $\R_+$.}

\vspace{0.2cm}\noindent
{\bf Theorem 3.}
 {\it Let $Y\sim G$ on $\R_+$ have all finite moments and  $G$ be absolutely continuous
on $\R_{+,x_0}$ with $G^{\prime}(x)=g(x)>0$ on $\R_{+,x_0}.$ Involve the functions $\phi$ and $\varphi$ as described in Assumption 1, and
assume further that
\begin{eqnarray}
	\gamma_3(g,\phi):=\limsup_{x\to\infty}\frac{g(x+\phi(x))}{g(x)}<1.
\end{eqnarray}
Then, both $Y$ and $Y^2$ satisfy Carleman's condition $({\rm C_S})$ and are M-det on $\R_+.$}

We provide now a corollary which extends Theorem 2.2 in  [10]. Our conclusion is about both $Y$ and $Y^2$, not only for $Y$.

\vspace{0.2cm}\noindent
{\bf Corollary 3.}
{\it Let $Y\sim G$ on $\R_+$ have all finite moments  and  $G$ be absolutely continuous
on $\R_{+,x_0}=[x_0, \infty)$ for some $x_0> 1$ with $G^{\prime}(x)=g(x)>0$ on $\R_{+,x_0}.$ Consider the function $\phi_1$ with $a>0,\ \alpha\in[0,1]$, see $(2)$,  and assume that
\begin{eqnarray*}
	\limsup_{x\to\infty}\frac{g(x+a\,(\log x)^{\alpha})}{g(x)}<1.
\end{eqnarray*}
Then, $Y$ and $Y^2$ each satisfy Carleman's condition $({\rm C_S})$;  hence, they both
are M-det on $\R_+$.}

By using the two other choices, (3) and (4), namely

\vspace{0.2cm}
\centerline{$\phi_2 (x) = a\,(\log x)^{\alpha}+\log\log x \quad \mbox{ or } \quad \phi_3 (x) = a\,(\log x)^{\alpha}/\log\log x, \quad x > \e,$}

\vspace{0.1cm}\noindent
we arrive at two corollaries formulated separately. In both,  $a>0,\ \alpha \in [0, 1].$

\vspace{0.2cm}\noindent
{\bf Corollary 4.}
 {\it Let $Y\sim G$ on $\R_+$ have finite moments  and  $G$ be absolutely continuous
on $\R_{+,x_0}=[x_0, \infty)$ for some $x_0> {\rm e}$ with $G^{\prime}(x)=g(x)>0$ on $\R_{+,x_0}.$  Then
\[
\limsup_{x\to\infty}\frac{g(x+a\,(\log x)^{\alpha}+\log\log x)}{g(x)}<1 \quad
\Rightarrow \quad \mbox{ both } \ Y \ \mbox{ and } \ Y^2  \ \mbox { are M-det on }  \R_+.
\]
}

\vspace{0.2cm}\noindent
{\bf Corollary 5.}
{\it Let $Y\sim G$ on $\R_+$ have all finite moments  and  $G$ be absolutely continuous
on $\R_{+,x_0}=[x_0, \infty)$ for some $x_0> {\rm e}$ with $G^{\prime}(x)=g(x)>0$ on $\R_{+,x_0}.$  Then
\[
\limsup_{x\to\infty}\frac{g(x+a\,(\log x)^{\alpha}/\log\log x)}{g(x)}<1
\quad
\Rightarrow \quad \mbox{ both } \ Y \ \mbox{ and } \ Y^2  \ \mbox { are M-det on }  \R_+.
\]
}

\section{Auxiliary Lemmas}\label{sec3}

To prove the main results stated in Section \ref{sec2}, we need two crucial lemmas which are slight modifications of,  e.g., Lemma 3.2 (pages 163 and 166), in [10].

\vspace{0.2cm}\noindent
{\bf Lemma 1.}
{\it For $n\in \N$ and real $\varepsilon>0,$  define the function $g(y)=n\log y-\varepsilon y,\ y>0.$ Then, for any integer
$n\ge (\varepsilon\,{\rm e})^{-1},$ we have two claims:\\
{\rm (i)} \ $\sup_{y>0}g(y)=n\log n-n\log (\varepsilon{\rm e})\le 2n\log n.$ \\
{\rm (ii)} \ If  $X\sim F$ on $\R$ has a finite moment of order $n$,  then for each $y_+>0,$\\
\[
n\int_{y_+}^{\infty}y^{n-1}(\log y){\rm d}F(y)\le 2n\log n\int_{y_+}^{\infty}y^{n-1}{\rm d}F(y)+
\varepsilon \int_{y_+}^{\infty}y^{n}{\rm d}F(y).
\]
}
{\bf Proof.}
Note that the function $g$ is twice differentiable on $(0,\infty)$ and
\[g^{\prime}(y)=n/y-\varepsilon,\ y>0, \quad  g^{\prime\prime}(y)=-n/y^2<0,\ y>0.\]
So,  $g$ is a concave function on $(0,\infty)$. Moreover, $g^{\prime}(n/\varepsilon)=0,$; hence, $g$ has a maximum value at $\bar{y}=n/\varepsilon.$
This proves Claim (i) from which Claim (ii) follows immediately.

\vspace{0.2cm}\noindent
{\bf Lemma 2.}
{\it Let $\{a_n\}_{n=1}^{\infty}$ be a positive sequence. If for some constants $b, c>0,$
\[
a_n\le c(n\log n)a_{n-1}+b^n,\ n=2,3,\ldots,
\]
then
\[
a_n\le d_0c^n(n\log n)^n,\ n=2,3,\ldots, \
\mbox{ where } \ d_0=a_1/c+\exp(b/c).
\]}
{\bf Proof.}
For simplicity, we use the symbol  `$!$'  as a `generalized factorial'
and denote

\smallskip
\centerline{$(n\log n)!=(n\log n)((n-1)\log (n-1))\cdots (2\log 2).$}

\smallskip
\noindent
Then by iteration, we obtain the following for $n\ge 2:$
\begin{eqnarray*}
a_n &\le & c(n\log n)a_{n-1}+b^n\le c(n\log n)[c((n-1)\log (n-1))a_{n-2}+b^{n-1}]+b^n\\
&= & c^2(n\log n)((n-1)\log (n-1))a_{n-2}+c(n\log n)b^{n-1}+b^n\\
&\le & c^{n-1}(n\log n)!\Big[a_1+\frac{c^{-1}b^2}{(2\log 2)!}+\frac{c^{-2}b^3}{(3\log 3)!}+\cdots+ \frac{c^{-(n-1)}b^n}{(n\log n)!}\Big]\\
&\le & c^n(n\log n)^n\Big[\frac{a_1}{c}+\frac{(b/c)^2}{2!}+\frac{(b/c)^3}{3!}+\cdots + \frac{(b/c)^n}{n!}\Big]\\
&\le & c^n(n\log n)^n[a_1/c+\exp(b/c)]=d_0c^n(n\log n)^n.
\end{eqnarray*}
The proof is complete.

\section{Proofs of the Main Results}\label{sec4}
\vspace{0.2cm}\noindent
{\bf Proof of Theorem 1.}
In parallel with the (usual) $n$th moment of $F$,  $m_n=\int_{-\infty}^{\infty}x^n{\rm d}F(x)$, we  now need its
 $n$th absolute moment denoted by $\mu_n=\int_{-\infty}^{\infty}|x|^n{\rm d}F(x)$. Furthermore, we split $\mu_n$ into two parts
and write $\mu_n=\mu_n^++\mu_n^-,$ where
\[\mu_n^+=\int_{[0,\infty)}|x|^n{\rm d}F(x)=\int_0^{\infty}x^n{\rm d}F(x), \quad
\mu_n^-=\int_{(-\infty,0)}|x|^n{\rm d}F(x)=\int_{-\infty}^0(-x)^n{\rm d}F(x).
\]
On the other hand, from condition (1), we have two inequalities:
\begin{eqnarray}\gamma_{1,+}(f,\phi):=\limsup_{x\to\infty}\frac{f(x+\phi(x))}{f(x)}\le \gamma_1(f,\phi)<1,
\end{eqnarray}
\begin{eqnarray}\gamma_{1,-}(f,\phi):=\limsup_{x\to -\infty}\frac{f(x-\phi(-x))}{f(x)}\le \gamma_1(f,\phi)<1.
\end{eqnarray}

We first sketch our proof as follows. The plan is to show, by condition (7), that
\begin{eqnarray}
\mu_n^+\le c_+(n\log n)\mu_{n-1}^++b_+^n,\ n=2,3,\ldots,
\end{eqnarray}
where $c_+$ and $b_+$ are positive constants independent of $n$. Similarly, by (8), we have
\begin{eqnarray}
\mu_n^-\le c_-(n\log n)\mu_{n-1}^-+b_-^n,\ n=2,3,\ldots,
\end{eqnarray}
where $c_-$ and $b_-$ are  positive constants independent of $n$.

Once this is down, combining (9) and (10) yields
\begin{eqnarray}
\mu_n=\mu_n^++\mu_n^-\le c(n\log n)\mu_{n-1}+b^n,\ n=2,3,\ldots,
\end{eqnarray}
where $c=\max\{c_+, c_-\}$ and $b=b_++b_-.$ Then, applying Lemma 2 to (11), we find that
\begin{eqnarray}
\mu_n\le d_0c^n(n\log n)^n,\ n=2,3,\ldots,
\end{eqnarray}
where $d_0=\mu_1/c+\exp(b/c).$ Consequently, from (12), it follows that
\begin{eqnarray}
\sum_{n=1}^{\infty}\mu_n^{-1/(2n)}\ge \sum_{n=2}^{\infty}\mu_n^{-1/(2n)}\ge \sum_{n=2}^{\infty}d_0^{-1/(2n)}c^{-1/2}(n\log n)^{-1/2}=\infty,
\end{eqnarray}
and also that
\begin{eqnarray}
\sum_{n=1}^{\infty}m_{2n}^{-1/(2n)}= \sum_{n=1}^{\infty}\mu_{2n}^{-1/(2n)}\ge \sum_{n=1}^{\infty}d_0^{-1/(2n)}c^{-1}(2n\log (2n))^{-1}=\infty.
\end{eqnarray}

From (13), we see that the random variable $|X|$ satisfies Carleman's condition $({\rm C_S})$ and
is M-det on $\R_+.$ Similarly, from (14), we conclude that
  the random variable $X$ satisfies Carleman's condition $({\rm C_H})$ and is M-det on $\R$, and
 $|X|^2=X^2$ satisfies Carleman's condition $({\rm C_S})$ and is M-det on $\R_+.$

Recall that $F$ is not a discrete distribution;
neither are $|X|$ and $X^2.$ Therefore, both $|X|$ and $X^2$ are also
M-det on $\R.$ This concludes the proof of the theorem.

Inequality (9) remains to be proven. By condition (7), we take $\beta=\frac12\,(1+\gamma_{1,+}(f,\phi)),$
and see that there exists a number  $\hat{x}_0\ge x_0$ such that
 two inequalities hold:
\[
\hat{y}_0=\hat{x}_0+\phi(\hat{x}_0)\ge y_*
\]
 (see the three conditions (a)--(c) in Assumption 1 for the point $y_*$), and
\[\frac{f(x+\phi(x))}{f(x)}\le\beta<1,\ \ x\ge \hat{x}_0.
\]
Therefore,
\[f(x)-f(x+\phi(x))\ge (1-\beta)f(x),\ \ x\ge \hat{x}_0.
\]
This implies that for $n\ge 2,$ we obtain a lower bound for the integral
$I(\hat{x}_0)$, where
\begin{eqnarray}
 I(\hat{x}_0)&:=& \int_{\hat{x}_0}^{\infty}x^n(f(x)-f(x+\phi(x))){\rm d}x\\
&\ge & (1-\beta)\int_{\hat{x}_0}^{\infty}x^nf(x){\rm d}x=(1-\beta)\mu_n^+-(1-\beta)\int_0^{\hat{x}_0}x^n{\rm d}F(x)\nonumber\\
&\ge & (1-\beta)\mu_n^+-\hat{x}_0^n.
\end{eqnarray}

On the other hand, to estimate the upper bound of  (15), we consider the integral
$J_n(\hat{y}_0)$, where
\begin{eqnarray}
J_n(\hat{y}_0)&:=&\int_{\hat{y}_0}^{\infty}[y^n-(y-\varphi(y))^n](1-\varphi^{\prime}(y))f(y){\rm d}y\nonumber\\
&\le & \int_{\hat{y}_0}^{\infty}ny^{n-1}\varphi(y)f(y){\rm d}y\\
&\le & C_+n\int_{\hat{y}_0}^{\infty}y^{n-1}(\log y)f(y){\rm d}y.
\end{eqnarray}
Here, inequality (17) follows from condition (a) and the    Mean-Value Theorem,  because
\[y^n-(y-\varphi(y))^n=\varphi(y)\, n \, \tilde{y}^{n-1}\le \varphi(y)\cdot ny^{n-1},
\]
where $\tilde{y}\in [y-\varphi(y), \ y],$
while the second inequality (18) follows from condition (b).

Now, we are ready to estimate the upper bound of the integral (15). By changing variables, we rewrite $I(\hat{x}_0)$ as follows:
\begin{eqnarray}
 I(\hat{x}_0)
&=& \int_{\hat{x}_0}^{\infty}x^nf(x){\rm d}x-\int_{\hat{x}_0}^{\infty}x^nf(x+\phi(x)){\rm d}x\nonumber\\
&=& \int_{\hat{x}_0}^{\infty}y^nf(y){\rm d}y-\int_{\hat{y}_0}^{\infty}(y-\varphi(y))^n(1-\varphi^{\prime}(y))f(y){\rm d}y\nonumber\\
&=& J_n(\hat{y}_0)+\int_{\hat{y}_0}^{\infty}y^n\varphi^{\prime}(y)f(y){\rm d}y+\int_{\hat{x}_0}^{\hat{y}_0}y^nf(y){\rm d}y\nonumber\\
&\le & C_+n\int_{\hat{y}_0}^{\infty}y^{n-1}(\log y)f(y){\rm d}y+C_+\int_{\hat{y}_0}^{\infty}y^{n-1}f(y){\rm d}y+\hat{y}_0^n\nonumber\\
&\le & C_+n\int_{\hat{y}_0}^{\infty}y^{n-1}(\log y)f(y){\rm d}y+C_+\mu_{n-1}^++\hat{y}_0^n,
\end{eqnarray}
where we applied inequality (18) and condition (c). Finally, applying Lemma 1(ii) to  inequality (19), we find that for real
$\varepsilon>0$ and integer
$n\ge (\varepsilon{\rm e})^{-1}\ge 2,$
\begin{eqnarray}
I(\hat{x}_0) &\le& 2C_+(n\log n)\int_{\hat{y}_0}^{\infty}y^{n-1}f(y){\rm d}y+
C_+\varepsilon \int_{\hat{y}_0}^{\infty}y^{n}f(y){\rm d}y+C_+\mu_{n-1}^++\hat{y}_0^n\nonumber\\
&\le & 2C_+(n\log n)\mu_{n-1}^++C_+\varepsilon \mu_n^++C_+\mu_{n-1}^++\hat{y}_0^n\nonumber\\
&\le & 3C_+(n\log n)\mu_{n-1}^++C_+\varepsilon \mu_n^++\hat{y}_0^n.
\end{eqnarray}

Combining inequalities (16) and (20), we find that for the previously chosen $\beta, \hat{x}_0, \hat{y}_0,$ and
$0< \varepsilon < \min \{(1-\beta)/C_+, 1/(2{\rm e})\},$ the following inequality holds:

\[
(1-\beta-C_+\varepsilon)\mu_n^+\le 3C_+(n\log n)\mu_{n-1}^++(\hat{x}_0^n+\hat{y}_0^n),\ \ n\ge (\varepsilon{\rm e})^{-1}\ge 2.
\]
Therefore,
\[
\mu_n^+\le \bar{c}(n\log n)\mu_{n-1}^++\bar{b}^n,\ \ n\ge (\varepsilon{\rm e})^{-1}\ge 2,
\]
where $\bar{c}=3C_+/(1-\beta-C_+\varepsilon)$ and $\bar{b}=(\hat{x}_0+\hat{y}_0)/(1-\beta-C_+\varepsilon).$ This means that there exists an integer  $n_0$, \ $2\le n_0\in\N,$ and
two positive constants
$\bar{b}, \ \bar{c}$ such that
\[\mu_n^+\le \bar{c}(n\log n)\mu_{n-1}^++\bar{b}^n,\ \ n\ge n_0.
\]
Putting the previous $n_0-1$ terms $\{\mu_1^+, \mu_2^+,\ldots,\mu_{n_0-1}^+\}$ together, we can choose (if necessary) two larger
 positive constants $b_+, c_+$ such that for all $n\ge 2$, we have
\[
\mu_n^+\le {c_+}(n\log n)\mu_{n-1}^++b_+^n.
\]
This proves inequality (9) and completes the proof of Theorem 1.

\vspace{0.2cm}\noindent
{\bf Proof of Corollary 1.} For $a>0,\ \alpha\in[0,1],$ consider the function $\phi(x)= \phi_1 (x) = a(\log x)^{\alpha}$ on $[x_0,\infty).$
Then $\phi (x)\ge 0$ is differentiable and $y(x)=x+\phi(x)$ is strictly increasing in $x$. Notice that the inverse function $x(y)$ exists. The functions $y(x)$ and
 $\varphi(y)=y-x(y)$ together obey the following properties:
 \[\lim_{x\to\infty}\frac{y(x)}{x}=1,\quad \lim_{x\to\infty}\frac{\log(y(x))}{\log x}=1, \quad \lim_{x\to\infty}\frac{(\log y(x))^{\alpha}}{(\log x)^{\alpha}}=1,
 \]
 \[\varphi^{\prime}(y)=\frac{{\rm d}\varphi(y)}{{\rm d}y}=\frac{{\rm d}\phi(x)}{{\rm d}x}\cdot\Big(\frac{{\rm d}y}{{\rm d}x}\Big)^{-1}=\frac{\phi^{\prime}(x)}{1+\phi^{\prime}(x)}
 \in[0,1],\quad y\ge y_0=y(x_0),
 \]
 \[\quad \lim_{y\to\infty}\frac{\varphi(y)}{(\log y)^{\alpha}}=a>0,\quad \lim_{y\to\infty}\frac{y\varphi^{\prime}(y)}{(\log y)^{\alpha-1}}=a\alpha>0.
 \]

 Consequently,
 \[
 \lim_{y \to \infty} \frac{\varphi (y)}{\log y} = a \quad  \mbox{ for } \ \alpha =1, \quad \mbox{ and } \ 0, \ \mbox{ if } \ \alpha < 1
 \]
 and
 \[
 \lim_{y \to \infty} y\,\varphi'(y) = a \quad  \mbox{ for } \ \alpha =1, \quad \mbox{ and } \ 0, \ \mbox{ if } \ \alpha < 1.
 \]

  Therefore, the function $\varphi$
satisfies conditions (a)--(c). This proves Corollary 1 for the  choice  $\phi = \phi_1$; see (2). Similar arguments are used for the two other choices, (3) and (4).

\vspace{0.2cm}\noindent
{\bf Proof of Theorem 2.}
We start with the distribution $G$ on $\R_+$ and
define the symmetric distribution $H$ on $\R$ as follows:
\[
H(x)=
\left\{
\begin{array}{cc}
\frac12(G(x^2)+1), & x\ge 0\\
-\frac12(G(x^2)-1), & x<0.
\end{array}
\right.
\]
Then, $H$ is absolutely continuous on the set $\R_{\sqrt{x_0}}=\R\setminus (-\sqrt{x_0},\sqrt{x_0})$ and its density is
\[
h(x)=H^{\prime}(x)=
\left\{
\begin{array}{cc}
g(x^2)x, & x\ge \sqrt{x_0}\\
-g(x^2)x, & \quad x \leq -\sqrt{x_0}.
\end{array}
\right.
\]
Clearly, of $H$ and $G$, each determines the other one. Consider a random variable  $Z \sim H$.
Then, recalling that $Y \sim G$, we find that for each $n\in\N,$
\[
\E[Z^{2n-1}]=0\quad {\rm and}\quad \E[Z^{2n}]=\E[Y^n].
\]

It follows from condition (5) that
\begin{eqnarray*}
\limsup_{|x|\to\infty}\frac{h(x+{\rm{sign}}(x)\phi(|x|))}{h(x)}&=&\limsup_{x\to\infty}\frac{h(x+\phi(x))}{h(x)}\\
&=&\limsup_{x\to\infty}\Big[\frac{g((x+\phi(x))^2)}{g(x^2)}\cdot\frac{x+\phi(x)}{x}\Big]\\
&= &\limsup_{x\to\infty}\frac{g((x+\phi(x))^2)}{g(x^2)}=\gamma_2(g,\phi)<1.
\end{eqnarray*}
By Theorem 1, $Z \sim H$ satisfies Carleman's condition $({\rm C_H})$, namely,
\[\sum_{n=1}^{\infty}(\E[Z^{2n}])^{-1/(2n)}=\infty.
\]
Equivalently,
\[\sum_{n=1}^{\infty}(\E[Y^{n}])^{-1/(2n)}=\sum_{n=1}^{\infty}(\E[Z^{2n}])^{-1/(2n)}=\infty.
\]
Therefore, $Y\sim G$ satisfies Carleman's condition $({\rm C_S})$ and is M-det on $\R_+.$ The proof is complete.

\vspace{0.2cm}\noindent
{\bf Proof of Corollary 2.}
The proof follows the same steps as those in  Corollary 1 and is omitted.

\vspace{0.2cm}\noindent
{\bf Proof of Theorem 3.}
In this Stieltjes case, denote $\nu_n :=\E[Y^n]$ for $n\in \N.$
In view of the proof of Theorem 1 (e.g., that condition (7) implies  inequality (9)), we know that condition (6) also implies the following inequality:
\begin{eqnarray}
\nu_n\le c(n\log n)\nu_{n-1}+b^n,\ n=2,3,\ldots,
\end{eqnarray}
where $c$ and $b$ are  positive constants independent of $n$. Then applying Lemma 2 to (21), we conclude that
\begin{eqnarray}
\nu_n\le d_0c^n(n\log n)^n,\ n=2,3,\ldots,
\end{eqnarray}
where $d_0=\nu_1/c+\exp(b/c).$ Consequently, it follows from (22), as before, that
\[\sum_{n=1}^{\infty}\nu_n^{-1/(2n)}=\infty \quad \mbox{ and } \quad  \sum_{n=1}^{\infty}\nu_{2n}^{-1/(2n)}=\infty.
\]
Thus, both $Y$ and $Y^2$ satisfy Carleman's condition $({\rm C_S})$ and are M-det on $\R_+.$ The proof is complete.

\vspace{0.2cm}\noindent
{\bf Proofs of Corollaries 3--5.}
The proofs are based on the steps  similar to those in Corollary 1 and are omitted.

\section{Comments and Illustrative Example}\label{sec5}

 Whether or not a distribution with all finite moments is M-det (uniquely determined) has been a profound question in mathematics for more than a century. However, the M-det property is not less important from the applied point of view, in which case the so-called `checkable conditions' are most useful. There are a variety of  sufficient or necessary conditions available in the literature for either M-det or for M-indet. The conditions are given in terms of the moments, of the
 densities, if they exist, of the distribution tails, or an appropriate combination of these.  Classical results and/or modern developments can be found
in several sources. Among them are the papers [8] and [9] and the book [6]; see also the references therein.

 For completeness of the general picture, it must be mentioned that  there are deep fundamental mathematical results of the sort `if  and only if';  however, due to their complexity, they fall into the group of `uncheckable conditions'. Readers interested in these aspects can consult the books
[1], [2] and [3].

The contribution of this paper is as follows: We present precise checkable conditions in terms of the probability density tails under which a distribution is M-det. The proofs are based on the classical Carleman's condition. In general, also seen in the example below, in order to claim the M-det property of a distribution,  we do not need to calculate the moments, etc. Our results allow
for concluding the M-det property of a distribution to check just one single asymptotic relation in terms of  the density `dropping speed'.

 \vspace{0.2cm}\noindent
{\bf Example.}
(a) \ {\it Consider a
 random variable $X$ obeying the Gumbel distribution, $X \sim F$, where
\[
F(x)=\Lambda(x)={\exp}\{-{\e}^{-x}\},\ x\in \R.
\]
Then, \ $X, \ X^2$ and $|X|$ are all M-det on $\R$.  }

\smallskip
\noindent
(b) \ {\it For a real number $c>0$, let us define the truncated random variable $X_c=X\cdot {\rm 1}_{[X\ge c]}.$
	Then each of \ $X_c$ and  \ $X_c^2$ is M-det on $\R_+$ and also on $\R.$
}

Case (a) is a Hamburger case and it can be easily shown that $X$ has finite moments of all orders $n\in \N$. In fact, the moment-generating function of $X$ exists and in terms of Euler's gamma function, we have

\vspace{0.2cm}
\centerline{$M(t) := {\E}[\e^{tX}] =$ $\Gamma (1-t)$, \ so \ $M(t) < \infty$ \ {\rm \ for \ all} \ $t<1.$}

Hence, $X\sim F$ is M-det on $\R.$ However, by using our Corollary 1, we can conclude 'more', namely that
all random variables $X, X^2$ and  $|X|$ are M-det on $\R.$ To see this, we need the  density function  $f = F' = \Lambda'$:
\[
f(x)={\e}^{-x}{\exp}\{-{\e}^{-x}\},\ x\in \R.
\]
We claim that $f$ satisfies Condition (D), i.e., the relation
\[\limsup_{|x|\to\infty}\frac{f(x+{\rm{sign}}(x)\phi(|x|))}{f(x)}=0<1,
\]
for any choice of the function $\phi$: $\phi(x)=\log x,$ or $\log x+\log\log x,$ or $(\log x)/\log\log x.$

The above relation can be easily checked for each choice of $\phi$. For example, if taking the simplest one, $\phi(x)=\log x,$ we have, on the right tail, that
\begin{eqnarray*}
	\limsup_{x\to\infty}\frac{f(x+\log x)}{f(x)}&=&\limsup_{x\to\infty}\frac{{\e}^{-(x+\log x)}{\exp}\{-{\e}^{-(x+\log x)}\}}
	{{\e}^{-x}{\exp}\{-{\e}^{-x}\}}\\
	&=& \limsup_{x\to\infty}{\e}^{-\log x}{\exp}\left\{ -{\e}^{-x}\left(x^{-1}-1\right)\right\}=0.
\end{eqnarray*}
Similarly, on the left tail,
\begin{eqnarray*}
	\limsup_{x\to -\infty}\frac{f(x-\log (-x))}{f(x)}&=&\limsup_{x\to -\infty}\frac{{\e}^{-(x-\log (-x))}{\exp}\{-{\e}^{-(x-\log (-x))}\}}
	{{\e}^{-x}{\exp}\{-{\e}^{-x}\}}\\
	&=& \limsup_{x\to -\infty}{\e}^{\log (-x)}
	{\exp}\{{\e}^{-x}\,(x+1)\}=0.
\end{eqnarray*}

Case (b) is a Stieltjes case.  Define first the number ${\bar p}_c$, the mass on the right tail of $X$ as follows:
\[
{\bar p}_c={\mathsf P}\{X\ge c\}=\int_c^{\infty}{\rm d} \Lambda(x).
\]
Then the `new' random variable $X_c=X\cdot {\rm 1}_{[X\ge c]}$  has a distribution,   $F_c$, with density function
\[
f_c(x)=\frac{1}{{\bar p}_c}{\e}^{-x}{\exp}\{-{\e}^{-x}\},\ x\ge c; \quad f_c (x) = 0, \ x < c.
\]
Again, $X_c$ has finite moments of all orders $n\in \N$. Using our Corollary 3, we conclude that both  random variables $X_c$ and $X_c^2$ satisfy Carleman's
condition $(\rm C_S)$ and hence are M-det on $\R_+$, because the density function $f_c$ of $X_c$ satisfies the relation:
\[
\gamma(f_c):=\limsup_{x\to\infty}\frac{f_c(x+\log x)}{f_c(x)}=\limsup_{x\to\infty}\frac{f(x+\log x)}{f(x)}=0<1.
\]
The same holds for the two other choices of $\phi$.

Notice that here we do not need to carry out all moments of $X_c$ and $X_c^2$ in order to check Carleman's
condition $(\rm C_S)$. We  make a definite conclusion about the M-det property  of  $X_c$ and $X_c^2$  simply by referring  to the inequality $\gamma(f_c)<1$ which requires \mbox{standard analysis.}

\vspace{6pt}

\vspace{0.2cm}\noindent
{\bf Acknowledgments.} The authors thank the three anonymous Referees for their positive feedback and especially for  giving a couple of relevant suggestions for achieving clearness. Compared with v1 (8 Jul 2025), in this version, v2 (30 Oct 2025), we grouped some of the corollaries, made the text more compact and added a comprehensive example showing how to use the results in this paper.

\end{document}